\documentclass[a4paper,12pt]{article}

\usepackage{amsfonts}
\usepackage{amssymb}
\usepackage{amsmath}
\usepackage{mathptmx}
\usepackage{color}
\usepackage{cancel}
\usepackage{bbm}
\usepackage{pstricks}
\usepackage{psfrag}
\usepackage{mathrsfs}

\def\ap{'\thinspace}

\newtheorem{lemma}{Lemma}[section]
\newtheorem{theorem}{Theorem}[section]
\newtheorem{remark}{Remark}[section]

\newtheorem{proposition}{Proposition}[section]
\newtheorem{definition}{Definition}[section]

\newtheorem{example}{Example}[section]
\def\be{\begin{equation}}
\def\ee{\end{equation}}
\def\bea{\begin{eqnarray}}
\def\eea{\end{eqnarray}}
\def\beann{\begin{eqnarray*}}
\def\eeann{\end{eqnarray*}}
\def\bsea{\begin{subeqnarray}}
\def\esea{\end{subeqnarray}}
\def\bmat{\left[ \begin{array}}
\def\emat{\end{array} \right]} 
\def\bsmat{\left[ \begin{smallmatrix}}
\def\esmat{\end{smallmatrix} \right]} 

\def\ap{'\thinspace}

\def\proof{\noindent{\bf{\em Proof:}\ \ }}
\def\QED{\mbox{\rule[0pt]{1.5ex}{1.5ex}}}
\def\endproof{\hspace*{\fill}~\QED\par\endtrivlist\unskip}
\def\endex{\hspace*{\fill}~$\square$\par\endtrivlist\unskip}
\newcommand{\real}{{\mathbb{R}}}
\newcommand{\complex}{{\mathbb{C}}}

\newcommand{\tp}{{^\tra}}

\def\gD{{\cal D}}

\def\gR{{\cal R}}

\def\gU{{\cal U}}
\def\gV{{\cal V}}

\def\gX{{\cal X}}

\def\Ricc{{\bf{{R}}}}

\def\tilda{{\!\!\!\!\phantom{P}^\thicksim}}
\def\mtilda{{\!\!\!\!\phantom{P}^{-\thicksim}}}

\newcommand{\ima}{\operatorname{im}}
\newcommand{\rank}{\operatorname{rank}}
\newcommand{\normrank}{\operatorname{normrank}}
\newcommand{\diag}{\operatorname{diag}}
\newcommand{\defi}{\stackrel{\text{\tiny def}}{=}}

\def\tra{{\scalebox{.6}{\thinspace\mbox{T}}}}

\definecolor{Royalblue}{cmyk}{1,0.30,0.2,0.2}

\begin{document}
\begin{titlepage}
\title{\vspace{15mm}
The Generalised Discrete Algebraic Riccati Equation \\
in LQ optimal control \thanks{Partially supported by the Italian Ministry for Education and Research (MIUR) under PRIN grant n. 20085FFJ2Z).}\vspace{10mm}}
\author{{\large Augusto Ferrante$^\ddagger$  \quad Lorenzo Ntogramatzidis$^\star$ }\\
       {\small $^\ddagger$Dipartimento di Ingegneria dell\ap Informazione}\\
        {\small    Universit\`a di Padova, via Gradenigo, 6/B -- 35131 Padova, Italy}\\
       {\small     {\tt augusto@dei.unipd.it}} \\ 
       {\small     $^\star$Department of Mathematics and Statistics}\\
       {\small     Curtin University, Perth WA, Australia.}\\
       {\small    {\tt L.Ntogramatzidis@curtin.edu.au}}
 }%
\maketitle
\begin{center}
\begin{minipage}{14.2cm}
\begin{center}
\bf Abstract
\end{center}
This paper investigates the properties of the solutions of the generalised discrete algebraic Riccati equation arising from the solution of the classic infinite-horizon linear quadratic control problem. In particular, a geometric analysis is used to study the relationship existing between the solutions of the generalised Riccati equation and the output-nulling subspaces of the underlying system and the corresponding reachability subspaces.
This analysis reveals the presence of a subspace that plays an important role in the solution of the related optimal control problem,
which is reflected in the generalised eigenstructure of the corresponding extended symplectic pencil.


In establishing the main results of this paper, several ancillay problems on the discrete Lyapunov equation and spectral factorisation are also addressed and solved.
\end{minipage}
\end{center}
\begin{center}
\begin{minipage}{14.2cm}
\vspace{2mm}
\end{minipage}
\end{center}
\thispagestyle{empty}
\end{titlepage}
\section{Introduction}
\label{secintro}
Ever since in the early sixties Kalman described in his pioneering papers \cite{Kalman-60-BSMM,Kalman-64} the crucial role of Riccati equations in the solution of the linear quadratic (LQ) optimal control and filtering problems, the range of control and estimation problems where Riccati equations have been discovered to play a fundamental role has been increasing dramatically. 
Indeed, in the last fifty years Riccati equations have been found to arise also in linear dynamic games with quadratic cost criteria, spectral factorisation problems, singular perturbation theory, stochastic realization theory and identification, boundary value problems for ordinary differential equations, invariant embedding and scattering theory. 
For this reason, Riccati equations are universally regarded as a cornerstone of modern control theory.
Several monographs have been entirely devoted to providing a general and systematic framework for the study of Riccati equations, see e.g. \cite{Willems-BL-91,Lancaster-95,Ionescu-OW-99,Abou-Kandil-FIJ-03}.

In the continuous time, the structure of the solution of a linear-quadratic problem strongly depends on the rank of the matrix penalising the control in the performance index, which is traditionally denoted by $R$. When $R$ is non-singular the optimal control can be found by solving a Riccati equation (which is differential or algebraic depending on the horizon of the performance index). 
Indeed, in this case, such Riccati equation -- which explicitly involves the inverse of $R$ -- is well-defined. 
But when $R$ is singular, a solution of the problem is guaranteed to exist for all initial conditions only if the class of allowable controls is extended to include distributions  \cite{Hautus-S-83,Willems-KS-86,Prattichizzo-MN-04}, and the Riccati equation is not defined.  \\

In the discrete time, the classic solution of the infinite-horizon LQ problem  is traditionally expressed in terms of the solution $X$ of the Riccati equation
\bea
\label{dare}
X \! =  \! A^\tra X A \! - \! (A^\tra X B \! + \! S)(R \! + \! B^\tra X B)^{-1}(B^\tra X A \! + \! S^\tra) \! + \! Q,
\eea
where $A \in \real^{n \times n}$, $B \in \real^{n \times m}$, $Q\in \real^{n \times n}$, $S \in \real^{n \times m}$ and $R \in \real^{m \times m}$ are such that 
\be
\label{popov}
\Pi \defi \bmat{cc} Q & S \\ S^\tra & R \emat =\Pi^\tra \ge 0.
\ee
Matrix $\Pi$ is usually referred to as {\em Popov matrix}. The set of matrices $\Sigma=(A,B;Q,R,S)$ is often referred to as {\em Popov triple}, see e.g. \cite{Ionescu-OW-99}. Equation (\ref{dare}) is the so-called Discrete Riccati Algebraic Equation DARE($\Sigma$).
 Notice that now it is not the inverse of $R$ that explicitly appears in the Riccati equation but the inverse of the term $R+B^\tra\,X\,B$, which can be non-singular even when $R$ is singular. Nevertheless, even though the distinction between the cases in which $R$ is invertible or singular needs not be considered, very often even in the discrete time it is assumed that $R$ is non-singular because this assumption considerably simplifies several underlying mathematical derivations. 

However, even the solution to the infinite-horizon LQ problem expressed in terms of matrices satisfying this equation is somehow restrictive. Indeed, 
an LQ problem may have solutions even if DARE has no solutions, and the optimal control can be written in this case as a state feedback written in terms of a matrix $X$ such that $R+B^\tra\,X\,B$ is singular and satisfies the more general Riccati equation
\bea
X  &  =  &   A^\tra X A \! - \! (A^\tra X B \! + \! S)(R \! + \! B^\tra X B)^\dagger(B^\tra X A \! + \! S^\tra) \! + \! Q, \label{gdare} \\
 &   &  \qquad \ker (R+B^\tra\,X\,B) \subseteq \ker (A^\tra\,X\,B+S), \label{kercond}
\eea
where the matrix inverse in DARE($\Sigma$) has been replaced by the Moore-Penrose pseudo-inverse, see \cite{Rappaport-S-71}. Equation (\ref{gdare}) is known in the literature as the {\em generalised discrete-time algebraic Riccati equation} GDARE($\Sigma$). The GDARE($\Sigma$) with the additional constraint on its solutions given by (\ref{kercond}) is sometimes referred to as {\em constrained generalised discrete-time algebraic Riccati equation}, herein denoted by CGDARE($\Sigma$). 
 It is obvious that (\ref{gdare}) constitutes a generalisation of the classic DARE($\Sigma$), in the sense that any solution of DARE($\Sigma$) is also a solution of GDARE($\Sigma$) -- and therefore also of CGDARE($\Sigma$) because (\ref{kercond}) is automatically satisfied since $\ker (R+B^\tra\,X\,B)=0_m$ -- but the {\em vice-versa} is not true in general. Despite its generality, this type of Riccati equation has not yet received a great deal of attention in the literature. It has only been marginally studied in the monographs \cite{Saberi-SC-95,Ionescu-OW-99,Abou-Kandil-FIJ-03} and in the paper \cite{Ferrante-04}. The only comprehensive contributions entirely devoted to the study of the solutions of this equation are \cite{Ionescu-O-96} and \cite{Stoorvogel-S-98}. The former investigates conditions under which the GDARE($\Sigma$) admits a stabilising solution in terms of the deflating subspaces of the extended symplectic pencil. The latter studies the connection between the solutions of this equation and the rank-minimising solutions of the so-called Riccati linear matrix inequality. In pursuing this task, the authors of \cite{Stoorvogel-S-98} derived a series of results that shed some light into the structural properties of the solutions of the generalised Riccati equation, and in particular in the fundamental role played by the term $R+B^\tra\,X\,B$. An example is the important observation according to which the inertia of this matrix $R+B^\tra\,X\,B$ -- that from now on we will denote by $R_X$ for the sake of conciseness -- is independent of the particular solution $X$ satisfying CGDARE($\Sigma$), \cite[Theorem 2.4]{Stoorvogel-S-98}. This implies that a given CGDARE($\Sigma$) cannot have one solution $X=X^\tra$ such that $R+B^\tra X\,B$ is non-singular and another solution $Y=Y^\tra$ for which $R+B^\tra Y\,B$ is singular. As such, {\bf i)} if $X$ is a solution of DARE($\Sigma$), then all solutions of CGDARE($\Sigma$) will also satisfy DARE($\Sigma$) and, {\bf ii)} if $X$ is a solution of 
CGDARE($\Sigma$) such that $R+B^\tra\,X\,B$ is singular, then DARE($\Sigma$) does not admit solutions. The results presented in \cite{Stoorvogel-S-98} are established in the very general setting in which the Popov matrix $\Pi$ is not necessarily positive semidefinite as in (\ref{popov}). 

In this paper we are interested in the connection of the use of the CGDARE($\Sigma$) in the solution of optimal control or filtering problems -- the so-called $H_2$-DARE in the terminology of \cite{Stoorvogel-S-98}. It is often taken for granted that the generalised discrete-time Riccati equation generalises the standard DARE($\Sigma$) in the solution of the infinite LQ optimal control problem in the same way in which \cite{Rappaport-S-71} established that the generalised Riccati difference equation generalises the standard Riccati difference equation in the solution of the finite-horizon LQ problem. However, to the best of the authors' knowledge, this fact has never been presented in a direct, self-contained and rigorous way. Thus, the first aim of this paper is to fill this gap, by showing in an elementary, yet rigorous, way, the connection of the CGDARE($\Sigma$) and the solution of the standard infinite-horizon LQ optimal control problem. The second aim of this paper is to provide a geometric picture describing the structure of the solutions of the CGDARE($\Sigma$) in terms of the output nulling subspaces of the original system $\Sigma$ and the corresponding reachability subspaces. 
Indeed, under the  usual assumption of positive semidefiniteness of the Popov matrix, the null-space of $R_X$ is independent of the solution $X$ of CGDARE($\Sigma$). Even more importantly, this null-space is linked to the presence of a subspace -- that will be identified and characterised in this paper -- which plays an important role in the characterisation of the solutions of CGDARE($\Sigma$), and also in the solution of the related optimal control problem. This subspace does not depend on the particular solution $X$, nor does the closed-loop matrix restricted to this subspace. This new geometric analysis will reveal that the spectrum of the closed-loop system is divided into two parts: the first depends on the solution $X$ of the CGDARE($\Sigma$), while the second -- coinciding exactly with the eigenvalues of the closed-loop restricted to this subspace -- is independent of it and does not appear in the generalised eigenstructure of the extended symplectic pencil. At first sight, this fact seems to constitute a limitation in the design of the optimal feedback, because it means that regardless of the solution of the generalised Riccati equation chosen for the implementation of the optimal feedback, the closed-loop matrix will always present a certain fixed eigenstructure as part of its spectrum. 
However, when $R+B^\tra\,X\,B$ is singular, the set of optimal controls presents a further degree of freedom -- which is also identified in \cite[Remark 4.2.3]{Saberi-SC-95} -- that allows to place all the poles of the closed-loop system at the desired locations without changing the cost. 

Several other important ancillary results of independent theoretical interest are derived in this paper. These include interesting considerations on the solutions of Hermitian Stein equations and spectral factorisation results that generalise the classic ones in more than one direction.

\section{Linear Quadratic optimal control and CGDARE}
In this section we analyse the connections between Linear Quadratic (LQ) optimal control and CGDARE.
Most of the results presented in this section are considered ``common wisdom". However, we have not been able to find a place where they have been derived in detail, so we believe that this section may be useful.
Consider the classic LQ optimal control problem. In particular, consider the discrete linear time-invariant system governed by
\be
\label{eqsys}
x_{t+1} = A\,x_t+B\,u_t,
\ee
where $A \in \real^{n \times n}$ and $B \in \real^{n \times m}$, and let the initial state $x_0\in \real^n$ be given. The problem is to find a sequence of inputs $u_t$, with $t = 0,1, \ldots,\infty$, minimising the cost function
\be
\label{cost}
J(x_0,u) \defi \sum_{t=0}^\infty \bmat{cc} x_t^\tra & u_t^\tra \emat \bmat{cc} Q & S \\ S^\tra & R \emat \bmat{c} x_t \\ u_t \emat.
\ee
 
 Before we introduce the solution of the optimal control problem, 
we recall some classic linear algebra results which will be useful in the sequel. We also give a proof of these results for the sake of completeness.

  \begin{lemma}
  \label{lem1}
  Consider the symmetric positive semidefinite matrix $P=\left[ \begin{smallmatrix} P_{11} & P_{12} \\[1mm] P_{12}^\tra & P_{22} \end{smallmatrix} \right]$. Then, 
  \begin{description}
  \item{\em {\bf (i)}} $\;\;\,\,\ker P_{12} \supseteq \ker P_{22}$;
  \item{\em {\bf (ii)}} $\;\,\,\,P_{12}\,P_{22}^\dagger\,P_{22}= P_{12}$;
  \item{\em {\bf (iii)}} $\;\,P_{12}\,(I-P_{22}^\dagger P_{22})=0$;
  \item{\em {\bf (iv)}} $\;\,P_{11}-P_{12} P_{22}^\dagger P_{12}^\tra \ge 0$.
  \end{description}
  \end{lemma}
  \proof
{\bf (i)} Since $P=P^\tra\ge0$, two matrices $C$ and 
$D$ exist such that $P=[\, C \;\;\; D \,]^\tra [\, C \;\;\; D \,]$ so that
$P_{12}=C^{\tra}\,D$ and $P_{22}\,{=}\,D^{\tra}\,D$. Let $x\,{\in}\,\ker\,P_{22}$. Then,
$0=x^\tra\,D^\tra D\,x=|| D\,x ||^2$, which gives $D\,x\,{=}\,0$. This in turn implies that $x\,{\in}\,\ker\,P_{12}$. 
{\bf (ii)} The inclusion $\ker P_{12} \supseteq \ker P_{22}$ can be rewritten as  $\ima\,P_{12}^{\tra}\,{\subseteq}\,\ima\,P_{22}$. Thus, a matrix $K\,{\in}\,\mathbb{R}^{n \times m}$ exists such that $P_{12}\,{=}\,K\,P_{22}$. On post-multiplying both sides of this identity by $P_{22}^\dagger\,P_{22}$ we obtain $P_{12}\,P_{22}^\dagger\, P_{22}=K\,P_{22}\,P_{22}^\dagger\,P_{22}=K\,P_{22}=P_{12}$.
{\bf (iii)}  Since as already proved $\ker P_{22} \subseteq \ker P_{12}$, a matrix $K$ exists such that $P_{12}=K\,P_{22}$. Therefore,
  $P_{12}\,(I - P_{22}^\dagger \, P_{22})=K\,P_{22}\,(I - P_{22}^\dagger \, P_{22})=K\,(P_{22} -P_{22}\, P_{22}^\dagger \, P_{22})=K\,(P_{22} -P_{22})=0$.
 {\bf (iv)} It follows directly from 
  $P_{11}-P_{12}\,P_{22}^\dagger P_{12}^\tra = 
  \bsmat I & -P_{12} \,P_{22}^\dagger \esmat \left[ \begin{smallmatrix} P_{11} & P_{12} \\[1mm] P_{12}^\tra & P_{22} \end{smallmatrix} \right] \bsmat I \\[1mm] -P_{22}^\dagger \,P_{12}^\tra \esmat \ge 0$.
\endproof 
  
  We now introduce some notation that will be used throughout the paper. First, to any matrix $X=X^\tra \in \real^{n \times n}$  we associate the following matrices:
\bea
Q_X &  \defi  &  Q \! + \! A^\tra X A \! - \! X, \quad \!\!\!\!\!
S_X   \defi   A^\tra X\, B \! + \! S,\quad\!\!\!\!\!
R_X   \defi   R \! + \! B^\tra X B, \\
\label{defgx}
G_X  &  \defi  &  I_m-(R+B^\tra\,X\,B)^\dagger (R+B^\tra\,X\,B)=I_m-R_X^\dagger R_X, \\
K_X  &  \defi  &  (R+B^\tra\,X\,B)^\dagger (B^\tra\,X\,A+S^\tra)=R_X^\dagger S_X^\tra, \label{KX} \\
A_X  &  \defi  &  A-B (R+B^\tra\,X\,B)^\dagger (B^\tra\,X\,A+S^\tra)=A-B\,K_X, \\
\Pi_X  &  \defi  &  \left[ \begin{array}{cc} Q_X & S_X \\ S_X^\tra & R_X \end{array} \right].
\eea 
The term $R_X^\dagger R_X$ is the orthogonal projector that projects onto $\ima R_X^\dagger=\ima R_X$ so that $G_X$ is  the orthogonal projector that projects onto $\ker R_X$. Hence,
 $\ker R_X=\ima G_X$.
When $X$ is a solution of CGDARE($\Sigma$), then $K_X$ is the corresponding gain matrix, $A_X$ the associated closed-loop matrix, and $\Pi_X$  is the  so-called  {\em dissipation matrix}.
It is easy to see that all symmetric and positive semidefinite solutions of GDARE($\Sigma$) satisfy (\ref{kercond}), and are therefore solutions of CGDARE($\Sigma$). In fact, if $X$ is positive semidefinite, 
we find
\[
\left[ \begin{array}{cc} Q_X+X & S_X \\ S_X^\tra & R_X \end{array} \right]=\left[ \begin{array}{cc} A^\tra \\ B^\tra \end{array} \right] X \left[ \begin{array}{cc} A & B \end{array} \right] +\Pi \ge 0. 
\]
Therefore, applying Lemma \ref{lem1} we find (\ref{kercond}), that can be rewritten as $\ker R_X \subseteq \ker S_X$ and also as $S_X\,G_X=0$.  

The following fundamental result holds.
 
 \begin{theorem}
 Suppose that for every $x_0$ there exists an input $u_t \in \real^m$, with $t \in \mathbb{N}$, such that
$J(x_0,u)$ is finite. Then we have:
\begin{enumerate} 
 \item CGDARE($\Sigma$) admits symmetric solutions: A solution $\bar{X}=\bar{X}^\tra\ge 0$  may be obtained as the limit of the sequence of matrices generated by iterating the {\em generalised Riccati difference equation} (see (\ref{defR}) below) with zero initial condition.
 \item The value of the optimal cost is {\em $x_0^\tra \bar{X} x_0$}.
 \item $\bar{X}$ is the minimum positive semidefinite solution  of CGDARE($\Sigma$).
 \item The set of {\em all} optimal controls minimising (\ref{cost}) can be parameterised as 
 \be
u_t=-K_{\bar X}\,x_t+G_{\bar X}\,v_t, \label{optcontr}
\ee
with arbitrary $v_t$.
\end{enumerate}
\end{theorem}
\proof 
(1). Consider the finite horizon LQ problem consisting in the minimisation of the performance index with zero terminal cost
\be\label{indexT}
J_T \defi  \sum_{t=0}^{T-1} \bmat{cc} x_t^\tra & u_t^\tra \emat \Pi \bmat{c} x_t \\ u_t \emat
\ee
subject to (\ref{eqsys}) with assigned initial state $x_0 \in \real^n$. The optimal control is obtained (see e.g. \cite{Rappaport-S-71}) by iterating, {backward in time starting from the terminal condition $P_T(T)=0$, the generalised Riccati difference equation $P_{T}(t)=\Ricc[ P_{T}(t+1) ]$, where $\Ricc[\cdot]$ is the Riccati operator defined as
\bea
\label{defR}
\Ricc[ P ] \defi
 A^\tra P A \! - \! (A^\tra P B \! + \! S)(R \! + \! B^\tra P B)^\dagger(B^\tra P A \! + \! S^\tra) \! + \! Q
\eea
and the optimal value of the cost is $J_T^*(x_0)=x_0^\tra P_T(0) x_0$. Let us  now consider the ``reverse time" sequence of matrices defined as $X_t\defi P_t(0)$. 
Since $P_\tau(t)=P_{\tau-t}(0)$ for all $t \le \tau$, the sequence $\{X_t\}_{t \in\mathbb{N}}$ is obtained by iterating the generalised Riccati difference equation forward with initial condition $X_0=0$. 
The sequence $\{J_t^*(x_0) \defi x_0^\tra X_t x_0\}_{t \in \mathbb{N}}$ is obviously monotonically non-decreasing (it is the sequence of optimal costs over intervals of increasing  lengths $t$). Hence, the sequences $\{X_t\}_{t \in \mathbb{N}}$, and  $\{R+B^\tra\,X_t\,B\}_{t \in \mathbb{N}}$
are monotonically non-decreasing sequences of positive semidefinite matrices.
We now show that these sequences are bounded.
Assume, by contradiction, that $\lim_{t\rightarrow +\infty}\|X_t\|=+\infty.$
The sequence $\{X_t^1=\frac{X_t}{\|X_t\|}\}_{t \in \mathbb{N}}$ is bounded. Thus, there exists a converging sub-sequence $\{X_{t_i}^1\}$. Let  $\bar{X}^1$ be its limit.
Clearly $\|\bar{X}^1\|=1$: let $x_0^1\in\mathbb{R}^n$ be such that $\|x_0^1\|=1$ and
$(x_0^1)\!\!^\tra \bar{X}^1 x_0^1=1$. 
Since we assumed that for any $x_0$ there exists a trajectory that renders $J$ defined in (\ref{cost}) finite, we have that
there exist a constant $m_0$ and an input trajectory $u^1$ such that
$J_{t_i}^*(x_0^1)\leq J(x_0^1,u^1)\leq m_0$, where the first inequality follows from the optimality of the cost  $J_{t_i}^*(x_0^1)$ and the fact that, for a given $u^1$, the index (\ref{cost}) is a sum of infinite non-negative terms which is greater than or equal to the sum of the first $t_i$ terms of the sum.
On the other hand we have 
$J_{t_i}^*(x_0^1)=\|X_{t_i}\| (x_0^1)\!^\tra X_{t_i}^1 x_0^1\rightarrow +\infty $, which is a contradiction.\\
Since $\{X_t\}_{t \in \mathbb{N}}$ is non-decreasing and bounded, it admits limit $\bar{X}$ for $t \to \infty$. Then, $\lim_{t \to \infty} X_t=\lim_{t \to \infty} X_{t+1}=\lim_{t \to \infty} \Ricc [X_{t}]=\bar{X}$. Thus, if $\lim_{t \to \infty} \Ricc [X_{t}]=
\Ricc [ \bar{X}]$, then $\Ricc [\bar{X}]=\bar{X}$, i.e. $\bar{X}$ is a positive semidefinite solution of CGDARE($\Sigma$). To prove that this is indeed the case, it is sufficient to show that $\lim_{t \to \infty} R_{X_{t}}^\dagger = R_{\bar{X}}^\dagger.$ In fact, the pseudo-inverse is the only possible source  of discontinuity in the Riccati iteration. To prove the latter equality, consider the sequence $\{R+B^\tra\,X_t\,B\}_{t \in \mathbb{N}}$.
Since it  is a monotonically non-decreasing sequence of positive semidefinite matrices, }the chain of inclusions
\beann
\ker (R \! + \! B^\tra X_0B) \supseteq \ker (R \! + \! B^\tra X_1 B)\supseteq \ker (R \! + \! B^\tra X_2 B)\supseteq \ldots
\eeann
holds. Clearly, there exist a $\bar{t}$ such that for any $t \ge \bar{t}$ this chain becomes stationary, i.e., for any $t \ge \bar{t}$ there holds $\ker (R+B^\tra\,X_t\,B)=\ker (R+B^\tra\,X_{t+1}\,B)$. This implies that a change of coordinates independent of $t$ exists such that {in the new basis
$R_{X_{t}}=R+B^\tra\,X_t\,B=\diag\{R^0_t,O\}$, where $\{R^0_t\}_{t \ge \bar{t}}$, is a non-decreasing sequence of positive definite  matrices.
Clearly, $\lim_{t \to \infty} R_{X_{t}} = R_{\bar{X}}$, so that, in this basis, $R_{\bar{X}}$ has the form
$R_{\bar{X}}=R+B^\tra\,\bar{X}\,B=\diag\{R^0,O\}$,
where $R^0 \defi \lim_{t \to \infty} R_t^0$. Moreover, since the sequence $\{R^0_t\}$  is non-decreasing, $R^0$ is also nonsingular, so that $(R^0_t)^{-1}\rightarrow (R^0)^{-1}.$
Thus, in the chosen basis we have indeed 
\[
R_{X_{t}}^\dagger=(R\!+\!B^\tra X_t B)^\dagger=\bmat{cc}\! (R^0_t)^{-1} & O \\ O & O \emat
\!\longrightarrow\!
 \bmat{cc} \!(R^0)^{-1} & O \\ O & O \emat=R_{\bar{X}}^\dagger.
\]
(2). Let 
\be
\label{infimu}
J^\circ(x_0) \defi \inf_{u} J(x_0,u).
\ee
Clearly, $J^\circ(x_0)\geq J_t^*(x_0)=x_0^\tra X_t\, x_0$. Then, by taking the limit, we get $J^\circ(x_0)\geq x_0^\tra \bar{X} x_0$.
We now show that the time-invariant feedback control $u^*_t \defi -K_{\bar{X}} x_t$ yields the cost $x_0^\tra \bar{X} x_0$, which is therefore the optimal value of the cost.
Consider the cost index $J_{T,\bar{X}} \defi J_T + x_T^\tra \bar{X} x_T$, where $J_T$ is defined in (\ref{indexT}).
It follows from \cite[Section II]{Rappaport-S-71}, see also \cite{kalman1964fundamental}, that an optimal control for this index is given by the time-invariant feedback $u^*_t=-K_{\bar{X}} x_t$ and the optimal cost does not depend on the length $T$ of the time interval and is given by 
$J_{T,\bar{X}}^*=x_0^\tra \bar{X} x_0$. Notice that for this conclusion we only need the fact that $\bar{X}$ is a positive semi-definite solution of CGDARE($\Sigma$).
Now we have
\bea
 x_0^\tra \bar{X} x_0 &\leq & J^\circ(x_0)\leq J(x_0,u^*)=\sum_{t=0}^\infty \bmat{cc} x_t^\tra & (u^*_t)^\tra \emat \Pi \bmat{c} x_t \\ u^*_t \emat    \nonumber   \\
&=&\lim_{T\rightarrow\infty} 
\sum_{t=0}^T \bmat{cc} x_t^\tra & (u^*_t)^\tra \emat \Pi \bmat{c} x_t \\ u^*_t \emat \nonumber \\
&=&
\lim_{T\rightarrow\infty} J_{T,\bar{X}}^* -x_T^\tra \bar{X} x_T\leq
\lim_{T\rightarrow\infty} x_0^\tra \bar{X} x_0 =x_0^\tra \bar{X} x_0. \label{chainofine-eq}
\eea
Comparing the first and last term of the latter expression we see that all the inequalities are indeed equalities, so that the infimum in (\ref{infimu}) is  a minimum and its value is indeed $x_0^\tra \bar{X} x_0$.

(3). 
Suppose by contradiction that there exist another positive semidefinite solution $\tilde{X}$ of CGDARE($\Sigma$) and a vector $x_0\in\real^n$ such that
$x_0^\tra\,\tilde{X}\,x_0 <  x_0^\tra\,\bar{X}\,x_0$. 
Take the time-invariant feedback $\tilde{u}_t=-K_{\tilde{X}} x_t$.
The same argument that led to (\ref{chainofine-eq}) now gives 
$J(x_0,\tilde{u})\leq x_0^\tra \tilde{X} x_0<  x_0^\tra\,\bar{X}\,x_0$, which is a contradiction because we have shown that $  x_0^\tra\,\bar{X}\,x_0$ is the optimal value of the cost function $J$.

(4).
Let ${\mathcal U}_0$ be the set of optimal control inputs  at time $t=0$. Let $u_0\in\real^m$ and $x_1=A\,x_0+B\,u_0$ be the corresponding state at $t=1$. {Clearly the optimal cost can  be written as}
\be\nonumber
J^{*}=x_0^\tra\, \bar{X}\,x_0=\bmat{cc} x_0^\tra & u_0^\tra \emat \bmat{cc}  \bar{X} & O \\ O & O \emat\bmat{c} x_0 \\ u_0 \emat.
\ee
{Moreover, $u_0\in{\mathcal U}_0$ if and only if the optimal cost can  be written in the following alternative form:}
\beann
J^* &=& x_1^\tra\, \bar{X}\,x_1+\bmat{cc} x_0^\tra & u_0^\tra \emat \Pi \bmat{c} x_0 \\ u_0 \emat \\
& =&
\bmat{cc} x_0^\tra & u_0^\tra \emat \left( \bmat{c} A^\tra \\ B^\tra \emat  \bar{X} \bmat{cc} A \,&\, B \emat + \Pi\right) \bmat{c} x_0 \\ u_0 \emat.
\eeann
By subtracting the first expression from the second, we get that $u_0\in{\mathcal U}_0$ if and only if
$
\bmat{cc} x_0^\tra & u_0^\tra \emat \bsmat Q_{\bar{X} } & S_{\bar{X} } \\[1mm] S_{\bar{X} }^\tra & R_{\bar{X} } \esmat \bsmat x_0 \\[1mm] u_0 \esmat=0.
$
Since $\bar{X}$, and hence $\Pi_{\bar{X}}=\bsmat Q_{\bar{X} } & S_{\bar{X} } \\[1mm] S_{\bar{X} }^\tra & R_{\bar{X} } \esmat$, are positive semidefinite, this is equivalent to
$\bsmat Q_{\bar{X} } & S_{\bar{X} } \\[1mm] S_{\bar{X} }^\tra & R_{\bar{X} } \esmat \bsmat x_0 \\[1mm] u_0 \esmat=0$. 
Finally, this is equivalent to $u_0=-R_{\bar{X} }^{\dagger} \,S_{\bar{X} }^\tra\,x_0+G_{\bar{X} }\,v_0$, where $v_0\in\real^m$ is arbitrary, because the columns of $G_{\bar{X} }$ form a basis for $\ker R_{\bar{X} }$.  By iterating this argument for all $t=1,2,\ldots$, we get (\ref{optcontr}).}
\endproof

%

\section{Preliminary technical results}
In this section, we present several technical results that will be used in the sequel. Most of these are ancillary results on the Stein equation and on spectral factorisation of independent interest.

\subsection{The Hermitian Stein equation}
In this section, we give some important results on the solutions $X$ of the so-called Hermitian Stein equation (known also as the discrete-time Lyapunov equation):
\be
\label{stein}
X=A^\tra X\,A+Q,
\ee
where $A,Q \in \real^{n \times n}$ and $Q=Q^\tra \ge 0$.

\begin{lemma}
\label{lemstein}
Let $X$ be a solution of the Hermitian Stein equation (\ref{stein}). Then, $\ker X$ is $A$-invariant and is contained in the null-space of $Q$.
\end{lemma}
\proof
Let $\lambda \in \complex$ be on the unit circle and such that $(A+\lambda\,I_n)$ is invertible. We can re-write (\ref{stein}) as
\be
\label{steina}
X=A^\tra X\,(A+\lambda\,I_n)-\lambda\,A^\tra\,X+Q,
\ee
 so that
\beann
(\lambda\,A^\tra+I_n)\,X  &  =  &  A^\tra X\,(A+\lambda\,I_n)+Q 
 = \lambda\,A^\tra X\,(\lambda^*\,A+I_n)+Q,
\eeann
since $\lambda$ is on the unit circle (which means in particular that $\lambda^*=\lambda^{-1}$). This is equivalent to
\bea
&& X\,(\lambda^*\,A+I_n)^{-1} \nonumber \\
&& = \lambda\,(\lambda A^\tra\!+\!I_n)^{-1} A^\tra X+(\lambda A^\tra\!+\!I_n)^{-1} Q\,(\lambda^* A\!+\!I_n)^{-1}.\label{hu1}
\eea
Let $\xi \in \ker X$. On pre-multiplying (\ref{hu1}) by $\xi^*$ and post-multiplying it by $\xi$, we obtain $\xi^* (\lambda\,A^\tra+I_n)^{-1}\,Q\,(\lambda^*\,A+I_n)^{-1} \xi=0$, and 
since  $(\lambda\,A^\tra+I_n)^{-1}\,Q\,(\lambda^*\,A+I_n)^{-1}$ is Hermitian and positive semidefinite, we get
\be
\label{view}
Q\,(\lambda^*\,A+I_n)^{-1} \,\xi=0.
\ee
Let us now post-multiply (\ref{hu1}) by the same vector $\xi$. We get $X\,(\lambda^*\,A+I_n)^{-1}\,\xi=0$, which means that $\ker X$ is $(\lambda^*\,A+I_n)^{-1}$-invariant. 
Hence, it is also $(\lambda^*\,A+I_n)$-invariant and therefore $A$-invariant.
In view of (\ref{view}), $\ker X=(\lambda^*\,A+I_n)^{-1}\ker X$ is also contained in the null-space of $Q$.
\endproof

We recall that equation (\ref{stein}) has a unique solution if and only if $A$ is {\em unmixed}, i.e. for all pairs $\lambda_1,\lambda_2\in\sigma(A)$ we have $\lambda_1\lambda_2\neq 1$. In this case we have the following result.
\begin{lemma}
Let $A$ be unmixed and $X$ be the unique solution of (\ref{stein}) where $Q=Q^\tra \ge 0$. Then, $\ker X$ is the unobservable subspace of the pair $(A,Q)$.
\end{lemma}
\proof
Let the pair $(A,Q)$ be written in the Kalman observability canonical form, i.e.,
$A=\bsmat A_{11} & O \\[1mm] A_{21} & A_{22} \esmat$ and $Q=\bsmat Q_{1} & O \\[1mm] O & O \esmat$, where the pair $(A_{11},Q_1)$ is completely observable. Let us write (\ref{stein}) in this basis. We find
\beann
\left[ \begin{array}{cc} X_{11} & X_{12} \\ X_{12}^\tra & X_{22} \end{array} \right]  &  =  &  
\left[ \begin{array}{cc} A_{11}^\tra & A_{21}^\tra \\ O &  A_{22}^\tra \end{array} \right]\!\!\left[ \begin{array}{cc} X_{11} & X_{12} \\ X_{12}^\tra & X_{22} \end{array} \right]\!\!\left[ \begin{array}{cc} A_{11} & O \\ A_{21} & A_{22} \end{array} \right]\!+\!\left[ \begin{array}{cc} Q_{1} & O \\ O & O \end{array} \right] \\
 &  =  &  
\left[ \begin{array}{cc} \star & \star \\ \star & A_{22}^\tra X_{22} A_{22} \end{array} \right]+\left[ \begin{array}{cc} Q_{1} & O \\ O & O \end{array} \right].
\eeann
Therefore, $X_{22}$ satisfies the homogeneous Hermitian Stein equation $X_{22}=A_{22}^\tra X_{22} A_{22}$. Since  $A$ is unmixed, the submatrix $A_{22}$ is unmixed, and therefore $X_{22}=0$ is the unique solution of $X_{22}=A_{22}^\tra X_{22} A_{22}$. As such, the Hermitian Stein equation in this basis can be simplified as
\be\nonumber
\left[ \begin{array}{cc} X_{11} & X_{12} \\ X_{12}^\tra & O \end{array} \right]  =  
\left[ \begin{array}{cc} \star & A_{11}^\tra X_{12}\,A_{22} \\ A_{22}^\tra X_{12}^\tra\,A_{11} &  O \end{array} \right]+\left[ \begin{array}{cc} Q_{1} & O \\ O & O \end{array} \right].
\ee
Again, since $A$ is unmixed, for all $\lambda_1\in {\sigma(A_{11})}$ and $\lambda_2\in {\sigma(A_{22})}$ we have that 
$\lambda_1\lambda_2\neq 1$, and the top-right block of the latter equation yields the unique solution $X_{12}=0$. Therefore, we get the following equation for $X_{11}$: $X_{11}=A_{11}^\tra X_{11}\, A_{11}+Q_1$. In view of Lemma \ref{lemstein} the unique solution $X_{11}$ of the latter equation
has trivial kernel because the pair $(A_{11},Q_1)$ is observable. This implies that $\ker X= \ima \left[ \begin{smallmatrix} O \\[1mm] I \end{smallmatrix} \right]$ where the partition is consistent with the block structure of $X$. On the other hand, this subspace is exactly the unobservable subspace of the pair $(A,Q)$.
\endproof

\begin{remark}
{\em When the solution $X$ of the Hermitian Stein equation (\ref{stein}) is not unique, by Lemma \ref{lemstein} the null-space of $X$ is still 
$A$-invariant and is contained in the null-space of $Q$, but it could be strictly contained into the unobservable subspace of the pair $(A,Q)$ -- which we recall is the largest $A$-invariant subspace contained in the null-space of $Q$ -- without necessarily being equal to it. Moreover, it is possible that none of the solutions of the Hermitian Stein equation are such that $\ker X$ coincides with the unobservable subspace of the pair $(A,Q)$. Consider for example the Hermitian Stein equation (\ref{stein}) with
$A=\bsmat 1 & 0 \\[1mm] 1 & 1 \esmat$ and $Q=\bsmat 1 & 0 \\[1mm] 0 & 0 \esmat$.
In this case, it is easy to see that the set of  all solutions of the Hermitian Stein equation (\ref{stein}) is 
$
\gX=\left\{ \bsmat \alpha & -\frac{1}{2} \\[1mm] -\frac{1}{2} & 0 \esmat\,: \;\; \alpha \in \real\right\}.
$
The null-space of any solution of the Hermitian Stein equation is zero. Hence, none of the solutions of the Hermitian Stein equation is such that its kernel is equal to $\ima \left[ \begin{smallmatrix} 0 \\[1mm] 1 \end{smallmatrix} \right]$.
}
\end{remark}

The last result we need is the following.
  
\begin{lemma}
\label{genius}
Let $A\in\real^{n\times n}$, $F\in\real^{n\times n}$, $B\in\real^{n\times m}$ and assume $X\in\real^{n\times n}$ is such that
{\em
\be
\label{equno}
\bmat{c} A^\tra \\ B^\tra \emat X\,F=\bmat{c} X \\  O \emat.
\ee
}
Then, {\em $B^\tra\,(A^\tra)^k\,X=0$} for all $k \ge 0$, i.e., $\ima X$ is contained in the unobservable subspace of the pair {\em $(A^\tra,B^\tra)$}.
\end{lemma}
\proof
We first prove that $B^\tra\,X=0$. Let us choose a basis in which $F$ is written as
$F=\diag\{ N, F_I\}$, where $N$ is nilpotent and $F_I$ is invertible. Let us decompose $X$ accordingly, i.e., $X=\bmat{cc} X_1 & X_2 \emat$. It is very easy to see that $A^\tra X_1\,N=X_1$ implies $X_1=0$. In fact, by multiplying such equation by $A^\tra$ and $N$ to the left and to the right, respectively, we obtain $X_1=(A^\tra)^k \, X_1\,N^k$ for all $k\ge 0$.
By choosing $k$ to be greater than the nilpotency index of $N$, we get $X_1=(A^\tra)^k \, X_1\,N^k=0$. From (\ref{equno}) we also obtain $B^\tra X_2\,F_I=0$, which implies $B^\tra\,X_2=0$ since $F_I$ is invertible. Therefore, $B^\tra\,X=0$. \\
The same argument can be iterated to prove that $B^\tra\,(A^\tra)^k\,X=0$ for all $k \ge 0$. Indeed, by pre-multiplying the first of (\ref{equno}) by $A^\tra$ we get $A^\tra (A^\tra X)\,F=A^\tra X$. By pre-multiplying the same equation by $B^\tra$, we get $B^\tra (A^\tra X)F=0$ since we already proved that $B^\tra X$ is zero. Hence, we can write these two equations as
$\bsmat A^\tra \\[1mm] B^\tra \esmat (A^\tra X)\,F=\bsmat A^\tra X \\[1mm] O \esmat$
and re-apply the same argument used above to show that $B^\tra\,A^\tra\,X=0$, and so on.
\endproof

\subsection{Spectral Factorisation}
Since as aforementioned the Popov matrix $\Pi$ is assumed symmetric and positive semidefinite, we can consider a factorisation of the form
\be
\label{pifact}
\Pi=\left[ \begin{array}{cc} Q & S \\ S^\tra & R\end{array} \right]=\left[ \begin{array}{cc} C^\tra \\ D^\tra\end{array} \right]\left[ \begin{array}{cc} C \,&\, D \end{array} \right],
\ee
where $Q=C^\tra C$, $S=C^\tra D$ and $R=D^\tra D$.  Let us define the rational matrix
$W(z)   \defi  C\,(z\,I_n-A)^{-1}B+D$. 
The spectrum $\Phi(z) \defi W^\tilda(z)\,W(z)$ -- where $W^\tilda(z) \defi W^\tra(z^{-1})$ --  associated with the Popov triple $\Sigma$ can be written as
\beann
\Phi(z)=\left[ \begin{array}{cc} B^\tra (z^{-1}\,I_n-A^\tra)^{-1} & I_n \end{array} \right]
\left[ \begin{array}{cc} Q & S \\ S^\tra & R \end{array} \right]\,\left[ \begin{array}{cc} (z\,I_n-A)^{-1}\,B \\ I_n \end{array} \right],
\eeann
which is also referred to as {\em Popov function} associated with GDARE($\Sigma$), \cite{Ionescu-OW-99}. 
%
 The matrix inequality for an unknown matrix $X=X^\tra$ of the form $\Pi_X \ge 0$ is called the {\em discrete Riccati linear matrix inequality}, and is herein denoted by DRLMI($\Sigma$).
Let us also define
\be\nonumber
L(X) \defi \Pi_X-\Pi=\left[ \begin{array}{cc} A^\tra\ X\,A-X & A^\tra\,X\,B \\ B^\tra X\,A & B^\tra X\,B \end{array} \right].
\ee
Notice that $L(X)$ is a linear function of $X$.

\begin{lemma}{{\bf (\cite[p.322]{Stoorvogel-S-98}, see e.g. \cite{Colaneri-F-SCL} for a detailed proof).}
For any $X=X^\tra \in\real^{n\times n}$,} there holds 
\bea
\label{alpha}
\Phi(z)=\left[ \begin{array}{cc} B^\tra (z^{-1}\,I_n-A^\tra)^{-1} & I_n \end{array} \right]
\Pi_X\,\left[ \begin{array}{cc} (z\,I_n-A)^{-1}\,B \\ I_n \end{array} \right].
\eea
\end{lemma}
%

\begin{theorem}
\label{thespectral}
Let $r$ denote the normal rank of the spectrum $\Phi(z)$.\footnote{The normal rank of a rational matrix $M(z)$ is defined as 
$\normrank M(z) \defi \max_{z \in \complex} \rank M(z)$. The rank of $M(z)$ is equal to its normal rank for all but finitely many $z \in \complex$.}
 If $X$ is a solution of CGDARE($\Sigma$), the rank of $R_X$ is equal to $r$. If $X$ is a solution of DRLMI($\Sigma$), the rank of $R_X$ is at most equal to $r$.
\end{theorem}
\proof
Let us consider $X=X^\tra$ such that $\Pi_X \ge 0$. By Lemma \ref{lem1}, this means in particular that {\bf (i)} $R_X$ is positive semidefinite, {\bf (ii)} $\ker S_X \supseteq \ker R_X$, and {\bf (iii)} $Q_X-S_X\, R_X^\dagger S_X^\tra$ is positive semidefinite. Notice that {\bf (iii)} means that $X$ satisfies the discrete Riccati inequality
 \beann
&& \gD(X) \defi A^\tra\,X\,A \\
&& \;\;\; -(A^\tra\,X\,B+S)(R+B^\tra\,X\,B)^\dagger(B^\tra\,X\,A+S^\tra)+Q-X \ge 0.
 \eeann
 Therefore, we can write $\gD(X)=H_X^\tra H_X$ for some matrix $H_X$, which leads to the expression
 {
 \be
 \Pi_X=\left[ \begin{array}{c} S_X\,R_X^\dagger \\ I \end{array} \right]
 R_X \left[ \begin{array}{cc}  R_X^\dagger S_X^\tra & I \end{array} \right]+\bmat{c} H_X^\tra \\ O \emat \bmat{cc} H_X & O \emat. \label{hu}
\ee
By plugging (\ref{hu}) into (\ref{alpha}) we see that the spectrum $\Phi(z)=W^\tilda (z) W(z)$ can be written as
\beann
\Phi(z)&=&W_1^\tilda(z) W_1(z)+W_2^\tilda(z) W_2(z) \\
&=&\bmat{cc} W_1^\tilda(z) & W_2^\tilda(z) \emat \bmat{c} W_1(z) \\ W_2(z)\emat,
\eeann
where $W_1(z)$ is given by 
\beann
W_1(z) &=& R_X^{\frac{1}{2}}  \left[ \begin{array}{cc}  R_X^\dagger S_X^\tra & I \end{array} \right]
\bmat{c} (z\,I_n-A)^{-1} B \\ I_n \emat \\
& =& R_X^{\frac{1}{2}}(R_X^\dagger\,S_X^\tra\,(z\,I_n-A)^{-1}B+I_m),
\eeann
and $W_2(z)$ is given by
\beann
W_2(z)=\bmat{cc} H_X & O \emat \bmat{cc} (z\,I_n-A)^{-1} B \\ I_n \emat=H_X\,(z\,I_n-A)^{-1} B.
\eeann
Notice that $W_1(z)= R_X^{\frac{1}{2}}T_X(z)$, where $T_X(z) \defi R_X^\dagger\,S_X^\tra\,(z\,I_n-A)^{-1}B+I_m$ is square and invertible for all but finitely many $z \in \complex$. Its inverse can be written as 
$T^{-1}_X(z)=I_m-R_X^\dagger S_X^\tra (z\,I_n-A_X)^{-1}B$. Thus, the normal rank of 
\beann
T_X^{\mtilda}(z) \Phi(z) T_X^{-1}(z) \! = \bmat{cc} R_X^{\frac{1}{2}}\,\, & T_X^{\mtilda}(z)\,W_2^\tilda(z) \emat \!\!\! \bmat{c}  \!  \!R_X^{\frac{1}{2}}  \!  \! \\  \!  \! W_2(z)\,T_X^{-1}(z)  \! \! \emat
\eeann
is equal to  the normal rank $r$  of $\Phi(z)$. Then, the rank of $R_X^{\frac{1}{2}}$, which equals that of $R_X$, is not greater than $r$.
Now consider the case where $X=X^\tra$ is a solution of CGDARE($\Sigma$). In this case, the term $H_X$ in (\ref{hu}) is zero, and therefore so is the rational function $W_2(z)$. 
 As such, ${W}_1(z)$ is a square spectral factor of $\Phi(z)$, i.e., 
$
W^\tilda(z)\,W(z)={W}_1^\tilda(z)\,{W}_1(z).$
Moreover, $T_X^{\mtilda}(z) \Phi(z) T_X^{-1}(z)=R_X$, which implies that when $X=X^\tra$ is a solution of CGDARE($\Sigma$), the rank of $R_X$ is exactly $r$. }
\endproof

\section{Geometric properties of the solutions of GDARE}
 The first aim of this section is to show that, given a solution $X$ of GDARE($\Sigma$)
\begin{itemize}
\item the subspace $\ker X$ is an output-nulling subspace for the quadruple $(A,B,C,D)$, i.e.,
\bea
\label{on}
\bmat{c} A \\ C \emat \ker X \subseteq (\ker X \oplus 0_p)+\ima \bmat{c} B\\ D \emat;
\eea
\item the gain $K_X$ is such that $-K_X$ is a {\em friend} of $\ker X$, i.e.,
\bea
\label{friend}
\bmat{c} A-B\,K_X \\ C-D\,K_K \emat \ker X \subseteq \ker X \oplus 0_p.
\eea
\end{itemize}
{ In the case where  $X=X^\tra$ is the  solution of GDARE($\Sigma$) corresponding to the optimal cost, these properties are intuitive. 
Indeed, on the basis of the optimality and of the fact that the cost cannot be smaller than zero in view of the positivity of the index, it is not  not difficult to prove 
that the following stronger result holds.
\begin{proposition}
Let $X$ be the minimal positive semidefinite solution of GDARE($\Sigma$).
Then $\ker X$ is the {\em largest} output-nulling subspace of the quadruple $(A,B,C,D)$. 
Moreover, $-K_X$ is the corresponding friend.
\end{proposition}
}
\proof Let $x_0\in\ker X$. 
Since the corresponding optimal cost is $J=x_0^\tra X\,x_0=0$, the initial state $x_0$ must belong to the largest output-nulling subspace of the quadruple $(A,B,C,D)$. {\em Vice-versa}, if we take a vector $x_0$ of the largest output-nulling subspace $\gV^\star$ of the quadruple $(A,B,C,D)$, by definition it is possible to find a control $u_k$ ($k \ge 0$) such that the state trajectory lies on $\gV^\star$ by maintaining the output at zero. This means that the corresponding value of the cost is zero. Hence, $x_0^\tra X\,x_0=0$ implies $x_0 \in \ker X$. The fact that $-K_X$ is a friend of $\ker X$ follows straightforwardly from the fact that if the initial state of the system lies on $\ker X$ and we assume by contradiction that $(A-B\,K_{X})\,x_0 \notin \ker X$, then the corresponding trajectory is not optimal because it is associated with a strictly positive value of the cost. Moreover, since the optimal cost is zero, we must have $(C-D\,K_{X})\ker X=0_p$.
\endproof
These ideas can be easily generalised to prove (\ref{on}) and (\ref{friend}) for any positive semidefinite solution $X=X^\tra \ge0$ of GDARE($\Sigma$).
Our aim is to prove a deeper geometric result:  (\ref{on}) and (\ref{friend}) hold for any  symmetric solution  $X$ of GDARE($\Sigma$).


\begin{theorem}\label{th-inv+friend}
Let $X$ be a solution of GDARE($\Sigma$). 
Then, $\ker X$ is an output-nulling subspace of the quadruple $(A,B,C,D)$ and $-K_X$ is a friend of $\ker X$, i.e., (\ref{on}) and (\ref{friend}) hold.
\end{theorem}
\proof
{ Since $X$ is a solution of GDARE($\Sigma$), the identity
\be
\label{clgdare}
X=A_X^\tra X \,A_X+Q_{0X}
\ee
holds, where $Q_{0X}\defi\left[ \begin{smallmatrix} I_n & -S_X\,R_X^\dagger \end{smallmatrix} \right] \left[ \begin{smallmatrix} Q & S \\ S^\tra & R \end{smallmatrix} \right] \left[ \begin{smallmatrix} I_n \\ -R_X^\dagger\,S_X^\tra \end{smallmatrix} \right]\geq 0$.} In view of Lemma \ref{lemstein}, $\ker X$ is $A_X$-invariant and is contained in the null-space of
$Q_{0X}$. {By factorising $\Pi$ as in (\ref{pifact}), we get $Q_{0X}=C_X^\tra C_X$ where
\be
\label{defcx}
C_X \defi \left[ \begin{array}{cc} C\, &\, D \end{array} \right] \left[ \begin{array}{cc} I_n \\ -R_X^\dagger\,S_X^\tra \end{array} \right]=C-D\,R_X^\dagger S_X^\tra.
\ee
Hence, the subspace $\ker X$ is also contained in the null-space of $C_X $ so that
 $\ker X$ is output-nulling for the quadruple $(A,B,C,D)$ and $-K_X$ is a friend of $\ker X$.}
\endproof

\ \\

Our aim at this point is to provide a full characterisation of the reachable subspace on $\ker X$, for reasons that will become clear in the sequel. Indeed, we will show that this subspace plays a crucial role in the solution of the associated optimal control problem. 
We recall the following definition.
{\begin{definition}\label{defrs}
The reachable subspace $\gR^\star_{\gV}$ on an output-nulling subspace $\gV$ is the subspace of the points of $\gV$ that can be reached from the origin along trajectories contained on $\gV$ by at the same time maintaining the output at zero.
\end{definition}
}
We will show that the reachable subspace $\gR^\star_{\ker X}$ on $\ker X$, coincides with the classic reachable subspace from the origin of the pair $(A_X,B\,G_X)$.
In order to prove this fact, we first need to show some important additional results on the solutions of CGDARE($\Sigma$). In particular, we now focus our attention on the term $R_X=R+B^\tra\,X\,B$. 
We see immediately that when $X$ is positive semidefinite, the null-space of $R_X$ is given by the intersection of the null-space of $R$ with that of $X\,B$. This result, which is very intuitive and easy to prove for positive semidefinite solutions of CGDARE($\Sigma$),  indeed holds for any solution $X$. However, in this case the proof -- which is divided between Lemma \ref{lemmaR} and Lemma \ref{lemmaXB} presented below -- is much more involved, and requires the machinery constructed in the first part of the paper. 

{
\begin{lemma}
\label{lemmaR}
Let $X=X^\tra$ be a solution of CGDARE($\Sigma$), $C_X$ be defined by (\ref{defcx}) and
\be\label{r0}\gR_0 \defi \ima \bmat{ccccccc} B G_X & A_X B G_X & A_X^2 B G_X & \ldots & A_X^{n-1} B G_X \emat.
\ee
 Then,
\be
\label{condR}
\ker R_X \subseteq \ker R,\qquad {\rm and}\qquad \gR_0 \subseteq\ker C_X.
\ee
\end{lemma}
}
\proof
{Since  the columns of $G_X$ (defined in (\ref{defgx})) span $\ker R_X$, we need to show that $R\,G_X=0$.
Recall from the proof of Theorem \ref{thespectral} that when $X=X^\tra$ is a solution of CGDARE($\Sigma$), the spectrum $\Phi(z)$, can be written as
 $\Phi(z)=W^{\tilda}(z)\,W(z)=T^{\tilda}_{X}(z)R_X T_X(z)$ where $T_X(z)=R_X^\dagger\,S_X^\tra\,(z\,I_n-A)^{-1}B+I_m$ is square and invertible for all but finitely many $z \in \complex$.
Hence, we have $R_X=[W(z) T^{-1}_X(z)]^{\tilda}[W(z) T^{-1}_X(z)]$ so that $R_XG_X=0$ implies $W(z) T^{-1}_X(z)G_X\equiv 0$.
Now recall from the proof of Theorem \ref{thespectral} that 
$T^{-1}_X(z)=I_m-R_X^\dagger S_X^\tra (z\,I_n-A_X)^{-1}B$ so that we can compute}
\beann
W(z) T^{-1}_X(z)=\left( C (z I_n\!-\!A)^{-1}B\!+\!D  \right) \left( I_m\!-\!R_X^\dagger S_X^\tra (zI_n\!-\!A_X)^{-1}B \right).
\eeann
Consider the following term of the product:
\be\nonumber
H(z)=C\,(z\,I_n-A)^{-1}B\,R_X^\dagger S_X^\tra (z\,I_n-A_X)^{-1}B.
\ee
By noticing that $B\,R_X^\dagger S_X^\tra =A-A_X=(zI_n-A_X)-(zI_n-A)$, we obtain
\be\nonumber
H(z)=C\,(z\,I_n-A)^{-1}B-C\,(z\,I_n-A_X)^{-1}B.
\ee
Hence,
\beann
W(z) T^{-1}_X(z)  &  =  &  W(z)-D\,R_X^\dagger S_X^\tra (z\,I_n-A_X)^{-1}B -H(z) \\
  &  =  &  D+(C-D\,R_X^\dagger S_X^\tra)(z\,I_n-A_X)^{-1} B \\
  &  =  &   D+C_X(z\,I_n-A_X)^{-1} B.
 \eeann
 {Since $W(z) T^{-1}_X(z)G_X$ is identically zero, it must be zero also when $z \to \infty$. In particular, $D\,G_X=0$, so that $R\,G_X=0$, which yields the first of (\ref{condR}). 
From $W(z) T^{-1}_X(z)G_X\equiv 0$ we also get $C_X(z\,I-A_X)^{-1} BG_X\equiv 0$ so that the reachable subspace of the pair $(A_X,BG_X)$, i.e. (\ref{r0}), is contained in $\ker C_X$ so that also the second of  (\ref{condR}) holds.
}
 \endproof

 In Lemma \ref{lemmaR} we have shown that $\ker R_X \subseteq \ker R$. Since $R_X=R+B^\tra X\,B$, it also straightforwardly follows that $\ker R_X \subseteq \ker (B^\tra X\,B)$ for any solution $X$ of CGDARE($\Sigma$). However, a stronger result holds, which says that $\ker R_X \subseteq \ker (X\,B)$. { This is an obvious consequence of Lemma \ref{lemmaR} for any solution $X\geq 0$, while it is a quite surprising and deep geometric result in the general case.}

\begin{lemma}
\label{lemmaXB}
Let $X=X^\tra$ be a solution of CGDARE($\Sigma$). Then,
\be
\label{condXB}
\ker R_X \subseteq \ker (X\,B).
\ee
\end{lemma}
\proof
From Lemma \ref{lemmaR}, if $v\in\ker R_X$, then $v\in\ker R\cap\ker (B^\tra X\,B)$.
We can select a change of coordinates in the input space $\real^m$ induced by {the $m \times m$ orthogonal matrix $T_X=\bmat{cc} T_{1X} & T_{2X} \emat$ where $\ima T_{1X}=\ima R_X$ and $\ima T_{2X}=\ima G_{X}=\ker R_X$.} In this basis $R_X$ is block-diagonal, with the first block being non-singular and the second being zero. Since $\ker R \supseteq \ker R_X$ as proved in Lemma \ref{lemmaR}, matrix $R$ in this basis has the form $R=\left[ \begin{smallmatrix}R_1 & 0\\[1mm] 0&0\end{smallmatrix} \right]$. {In the same basis, matrix $B$ can} be partitioned accordingly as $B=\bmat{cc} B_1 & B_2\emat$, so that $\ima B_2=\ima (B\,G_X)$.
We must show that $X\,B_2=0$.
Since $\ker R_X \subseteq \ker (B^\tra\,X\,B)$, in this basis we find
\be
\label{bxb2=0}
B^\tra X B_2=\bmat{c}B_1^\tra\\ B_2^\tra\emat X\,B_2=0.
\ee
Moreover, since $\ker R \subseteq \ker S$, in the selected basis $S$ takes the form $S=\bmat{cc} S_1 & 0\emat$.
Thus, $S_X = A^\tra X\,B+S=\bmat{cc} A^\tra X B_1+S_1&A^\tra X B_2\emat$.
From $\ker R_X\subseteq \ker S_X$ it now follows that
$A^\tra X\,B_2=0$ which, together with (\ref{bxb2=0}), yields
\be
\label{AB}
\bmat{cc}
A^\tra\\ 
B^\tra
\emat X\,B_2=0.
\ee
If $A$ is non-singular or, more in general, if the zero eigenvalue of $A$, when present, is controllable from $B$, then clearly
$X\,B_2=0$. However, this result is true in general, without any assumption. To prove this, let us {consider $\gR_0$  defined in (\ref{r0}) which, in the chosen input space basis, is the reachable subspace of the pair $(A_X, B_2)$.\footnote{In the symbol denoting this subspace we dropped the subscript $X$ because, as it will be proved in the sequel, this subspace is independent of the particular solution of the CGDARE($\Sigma$).} 
Let us consider a basis of the state-space {where the pair $(A_X,B_2)$ are in Kalman controllability form. In such a basis,  the subspace $\gR_0$ is spanned by the columns of the matrix $\left[ \begin{smallmatrix} I \\[1mm] O  \end{smallmatrix} \right]$ and} we have
\bea
\label{dec}
A_X=\bmat{cc}  A_{X,11}   &   A_{X,12}   \\   O   &   A_{X,22}   \emat\!, \!\quad       B_2=\bmat{cc}     B_{21}     \\     O     \emat\!,\!       \quad B_1=\bmat{cc}     B_{11}     \\     B_{12}     \emat\!,
\eea
where the pair $(A_{X,11},B_{21})$ is reachable.
In this basis, matrix $C_X$ takes the form $C_X=[O \mid C_{X,1} ]$ in view of the second of (\ref{condR}).}
  Since $A_X=A-B\,K_X$, we can  re-write (\ref{AB}) as
$
\bsmat A_X^\tra +K_X^\tra B\tp \\[1mm]
B\tp \esmat X\,B_2=0
$ or, equivalently, as 
$\left[ \begin{smallmatrix} A_X^\tra \\[1mm] B^\tra \end{smallmatrix} \right]X\,B_2=0$. Using the partitioned structure described above, we can re-write this equation as
\be
\label{gamma1}
\bmat{cc} A_{X,11}^\tra & O \\ 
 A_{X,12}^\tra &  A_{X,22}^\tra \\
  B_{11}^\tra & B_{12}^\tra \\ B_{21}^\tra & O \emat \bmat{cc} X_{11} & X_{12} \\ X_{12}^\tra & X_{22} \emat \bmat{cc} B_{21} \\ O \emat=0.
  \ee
We want to show that $\bsmat X_{11} & X_{12} \\[1mm] X_{12}^\tra & X_{22} \esmat \bsmat B_{21} \\[1mm] O \esmat=0$, i.e., that $X_{11}\,B_{21}=0$ and $X_{12}^\tra B_{21}=0$. From (\ref{gamma1}) we find
  \bea
  A_{X,11}^\tra\,X_{11}\,B_{21}  &  =  &  0, \label{eq*} \\
  B_{21}^\tra\,X_{11}\,B_{21}  &  =  &  0. \label{eqo}
  \eea
Since the pair $(A_{X,11},B_{21})$ is reachable by construction, $X_{11}\,B_{21}=0$. It remains to show that $X_{12}^\tra B_{21}=0$.
In this basis, equation (\ref{clgdare}) -- which is exactly GDARE($\Sigma$) -- now reads as
\bea
\bmat{cc} X_{11} & X_{12} \\ X_{12}^\tra & X_{22} \emat &=&\bmat{cc} A_{X,11}^\tra & O \\ 
 A_{X,12}^\tra &  A_{X,22}^\tra \emat \bmat{cc} X_{11} & X_{12} \\ X_{12}^\tra & X_{22} \emat 
\bmat{cc} A_{X,11} & A_{X,12} \\ O &  A_{X,22} \emat \nonumber \\
&& +\bmat{c} O \\ C_{X,1}^\tra \emat \bmat{cc} O & C_{X,1} \emat,\label{GDARE2}
\eea
from which we find in particular $X_{11}=A_{X,11}^\tra X_{11} A_{X,11}$. This equation can be written together with (\ref{eq*}) and (\ref{eqo}) as 
\be\nonumber
\bmat{cc} 
A_{X,11}^\tra X_{11} A_{X,11}-X_{11} & A_{X,11}^\tra\,X_{11}\,B_{21} \\
B_{21}^\tra X_{11}\,A_{X,11} & B_{21}^\tra\,X_{11}\,B_{21} \emat =0.
\ee
Since the pair $(A_{X,11},B_{21})$ is reachable, we can apply Lemma 2.9 in \cite{Stoorvogel-S-98}, which guarantees that $X_{11}$ is zero. Now we can re-write (\ref{GDARE2}) as
\bea
\bmat{cc} O & X_{12} \\ X_{12}^\tra & X_{22} \emat& =&\bmat{cc} A_{X,11}^\tra & O \\ 
 A_{X,12}^\tra &  A_{X,22}^\tra \emat \bmat{cc} O& X_{12} \\ X_{12}^\tra & X_{22} \emat 
\bmat{cc} A_{X,11} & A_{X,12} \\ O &  A_{X,22} \emat \nonumber \\
&& +\bmat{cc} O & O \\ O & C_{X,1}^\tra C_{X,1} \emat.\label{GDARE3}
\eea
In particular, we get 
\be
X_{12}=A_{X,11}^\tra X_{12} A_{X,22}. \label{alpha1}
\ee
Moreover, by { plugging  the value $X_{11}=0$ into (\ref{gamma1}), after transposition, we obtain}  
\be
B_{21}^\tra X_{12}A_{X,22}=0. \label{beta1}
\ee
Equations (\ref{alpha1}) and (\ref{beta1}) can be re-written as
\be
\label{34}
\bmat{c} A_{X,11}^\tra \\ B_{21}^\tra \emat X_{12} A_{X,22}= \bmat{c} X_{12} \\ O \emat.
\ee
By applying the result in Lemma \ref{genius} for $k=0$, we get $B_{21}^\tra\,X_{12}=0$.
\endproof
\begin{remark}
{\em 
{In the last line of the proof of Lemma \ref{lemmaXB} we can applying Lemma \ref{genius} for $k \ge 0$, and obtain $X_{12}^\tra A_{X,11}^k B_{21}=0$ for all $k\ge 0$. Since the pair $(A_{X,11},B_{21})$ is reachable, this yields $X_{12}=0$. Therefore the following stronger result holds.}
}
\end{remark}

{\begin{proposition}
Let $X=X^\tra$ be a solution of CGDARE($\Sigma$) and 
$\gR_0$ be defined by (\ref{r0}).
 Then, $X\,\gR_0=0_n$.
\end{proposition}}

\begin{remark}
{\em 
{As an obvious corollary of Lemmas \ref{lemmaR} and \ref{lemmaXB}, we have that
\be
\label{ckrx}
\ker R_X =  \ker (XB) \cap \ker R=\ker  \left[ \begin{smallmatrix} X\,B \\[1mm] R \end{smallmatrix} \right].
\ee
}
}
\end{remark}
  
\begin{remark}
{\em The result established in Lemma \ref{lemmaXB} 
does not continue to hold if we only assume that $X=X\tp$ is a solution of 
 GDARE($\Sigma$) (instead of being a solution  of CGDARE($\Sigma$)). Consider for example the case
$A=\bsmat -1 & 0 \\[1mm] -5 & -6 \esmat$, $B=\bsmat -4 & 0 \\[1mm] 0 & -2 \esmat$, $C=\bsmat 0 & 1 \esmat$ and $D=\bsmat 4 & 0 \esmat$. 
It can be easily verified that $X=\diag \{-1,1\}$ is a solution of the GDARE($\Sigma$) but not of the CGDARE($\Sigma$). In this case, 
$
\ker R_X=\ima \left[ \begin{smallmatrix} 1 \\[1mm] 0 \end{smallmatrix} \right] \neq \ker  \left[ \begin{smallmatrix} X\,B \\[1mm] R \end{smallmatrix} \right]= 0_m.
$
}
\end{remark}

\begin{remark}
{\em The result established in Lemma \ref{lemmaXB} does not hold when the Popov matrix is not positive semidefinite. Consider the following numerical example in which
$A=\bsmat 1 & 1 \\[1mm] 1 & 1 \esmat$, $B=\bsmat 0 & 1 \\[1mm] 0 & 1 \esmat$, $Q=\bsmat 1 & 0 \\[1mm] 0 & -1  \esmat$, $S=\bsmat 0 & 0 \\[1mm] 0 & 0 \esmat$ and $R=\bsmat 1 & 0 \\[1mm] 0 & 0  \esmat$. A solution of CGDARE($\Sigma$) is given by $X=\diag\{1,-1\}$. Indeed, $X$ satisfies (\ref{kercond}) since one can easily verify that $S_X$ is the zero matrix and $\ker R_X=\ima \left[ \begin{smallmatrix} 0\\[1mm] 1 \end{smallmatrix} \right]$. Changing coordinates in the input space as shown in the proof of Lemma \ref{lemmaXB} leads to $B_1=\left[ \begin{smallmatrix} 0 \\[1mm] 0 \end{smallmatrix} \right]$ and $B_2=\left[ \begin{smallmatrix} 1 \\[1mm] 1 \end{smallmatrix} \right]$. However, this time $X\,B_2=\left[ \begin{smallmatrix} 1 \\[1mm] -1 \end{smallmatrix} \right]$.
}
\end{remark}

\begin{theorem}
\label{the}
Let $X=X^\tra$ be a solution of CGDARE($\Sigma$). Let $\gR_0$ denote the reachable subspace of the pair $(A_X,B\,G_X)$ as defined in (\ref{r0}), {and $\gR^\star_{\gV}$ be defined by Definition \ref{defrs}.} Then,
\be
\gR^\star_{\ker X}=\gR_0.
\ee
\end{theorem}
\proof
Let us first show that 
\be
\ima (B\,G_{X})=\ker X \cap B\,\ker D.
\ee
We recall that $\ima G_X=\ker R_X$. Moreover, from (\ref{ckrx})  we know that $\ker R_X=\ker (XB) \cap \ker R$.
Then $\ima (B\,G_{X})=B\ker R_X =  B(\ker (XB) \cap \ker R)= \ker X  \cap B\ker R  = \ker X  \cap B\ker D.$

Now we are ready to prove the statement of this theorem. 
Since $\gR_0$ is the reachable subspace from the origin of the pair $(A_X,B\,G_X)$, it is by definition the smallest $A_X$-invariant subspace containing $\ima (B\,G_X)=\ker X \cap B\,\ker D$. On the other hand, the reachable subspace $\gR^\star_{\ker X}$ on $\ker X$ is  {characterised as
follows \cite[Theorem  7.14]{Trentelman-SH-01}, \cite[p.~424]{Ntogramatzidis-05}: Let $F$ be an {\em arbitrary} friend of $\ker X$, i.e., $F$ is any feedback matrix such that 
$(A+B\,F)\ker X \subseteq  \ker X$ and $(C+D\,F)\ker X  =  0_p$.
Then $\gR^\star_{\ker X}$ is the smallest $(A+B\,F)$-invariant subspace 
containing $\ker X \cap B\,\ker D$. Notice that $\gR^\star_{\ker X}$ does not depend on the choice of the friend $F$, \cite[Theorem  7.18]{Trentelman-SH-01}.
We have seen in Theorem \ref{th-inv+friend} that the matrix $F= -K_X$ is a particular friend of $\ker X$.
For this choice of $F$, we have $A+B\,F=A-BK_X=A_X$, so that $\gR^\star_{\ker X}$ is the smallest $A_X$-invariant subspace 
containing $\ker X \cap B\,\ker D$, which is exactly the definition of $\gR_0$.}
\endproof

\begin{remark}
{\em The statement of Theorem \ref{the} can be also captured as follows. First, we know from Lemma \ref{lemmaR} that
$C_X(zI_n-A_X)^{-1}B\,G_X$ is identically zero, which implies that the reachable subspace from the origin of the pair $(A_X,B\,G_X)$ is contained in the non-observable subspace of the pair $(A_X,C_X)$. Thus, $\gR_0$ is also a controllability subspace for the quadruple $(A_X,0,C_X,0)$ and a controllability subspace for the quadruple $(A,B,C,D)$. Consider the same orthogonal change of coordinates in the input space introduced in the proof of Lemma \ref{lemmaXB} induced by the $m \times m$ matrix {$T_X=\bmat{cc} T_{1X}& T_{2X} \emat$.}  Let matrix $B$ be partitioned in this basis as $B=\bmat{cc} B_1 & B_2\emat$. This change of coordinates gives rise to a state-space model of the form
\beann
x_{k+1}  &  =  &  A_X\,x_k+\bmat{cc} B_1 & B_2 \emat v_k, \\
y_k  &  =  &  C_X\,x_k+\bmat{cc} D_1 & O \emat v_k. 
\eeann
By using the control $v_k=\left[ \begin{smallmatrix} 0 \\[1mm] \bar{v}_k \end{smallmatrix} \right]$, we find that the state $x_k$ can reach every point of the reachable subspace from the origin of the pair $(A_X,B_2)$ and at the same time the output $y_k$ is kept at zero. Therefore, this subspace is also a controllability subspace for $(A_X,0,C_X,0)$.
}
\end{remark}

 In \cite{Stoorvogel-S-98} it is proved that the inertia of $R_X$ is independent of the particular solution $X=X^\tra$ of CGDARE($\Sigma$). Here, we want to show that much more is true when $\Pi$ is positive semidefinite. Namely, the null-space of $R_X$ is independent of the particular solution $X=X^\tra$ of CGDARE($\Sigma$).
 
 \begin{theorem}
 \label{thker}
 Let $X_1,X_2$ be two solutions of CGDARE($\Sigma$). Then, $\ker R_{X_1}=\ker R_{X_2}$.
 \end{theorem}
 \proof
Consider two solutions $X_1=X_1^\tra$ and $X_2=X_2^\tra$ of CGDARE($\Sigma$). In particular, $X_1$ and $X_2$ also satisfy the generalised Riccati inequality, so that $\Pi_{X_1} \ge 0$ and $\Pi_{X_2} \ge 0$. In other words, using the same notation employed in the proof of Theorem \ref{thespectral}, we have $\Pi_{X_i}=L(X_i)+\Pi \ge 0$ for $i \in \{1,2\}$. The set of solutions of the generalised Riccati inequality is a convex set, i.e., by taking $\alpha \in (0,1)$, then $\alpha\,\Pi_{X_1}+(1-\alpha)\,\Pi_{X_2}$ is positive semidefinite because it is a convex combination of positive semidefinite terms. For our purposes, it is sufficient to fix an arbitrary value of $\alpha$, say $\alpha=\frac{1}{2}$. Then,
\beann
0 \le \frac{1}{2}\,(\Pi_{X_1}+\Pi_{X_2}) = \Pi_{\frac{1}{2}(X_1+X_2)}  =  \Pi+L\left(\frac{1}{2}(X_1+X_2)\right),
\eeann
where the last equality holds in view of the linearity of $L(\cdot)$. 
This means that $Y \defi \frac{1}{2}\,(X_1+X_2)$  satisfies the Riccati inequality $\Pi_Y \ge 0$. By virtue of Theorem \ref{thespectral}, the rank of $R_Y \defi R+B^\tra Y\,B$ is not greater than $r$. On the other hand,
\be
R_Y    =    R+\frac{1}{2}\,B^\tra (X_1+X_2) B = \frac{1}{2}\,(R_{X_1}+R_{X_2}). \label{eqfin}
\ee
Hence, since $X_1$ and $X_2$ are both solutions of CGDARE($\Sigma$), the ranks of $R_{X_1}$ and $R_{X_2}$ are exactly equal to $r$. Thus, the rank of $R_Y$ is greater or equal to $r$. This means that the rank of $R_Y$ must be exactly equal to $r$, i.e., from (\ref{eqfin}) we have that $R_{X_1}$ and $R_{X_2}$ must have the same null-space.
\endproof

Now we want to prove that the subspace $\gR^\star_{\ker X}$ is independent of the particular solution $X=X^\tra$ of CGDARE($\Sigma$). Moreover, $A_X$ restricted to this subspace does not depend on the particular solution $X=X^\tra$ of CGDARE($\Sigma$).

\begin{theorem}
Let $X$ and $Y$ be two solutions of CGDARE($\Sigma$). Let $A_X$ and $A_Y$ be the corresponding closed-loop matrices. Then,
\begin{itemize}
\item $\gR^\star_{\ker X}=\gR^\star_{\ker Y}$, and
\item {$ A_X|_{\gR^\star_{\ker X}}= A_Y|_{\gR^\star_{\ker Y}}$. } 
\end{itemize}
\end{theorem}
\proof
Let $\Delta \defi Y-X$. Since $\ker R_X$ coincides with $\ker R_Y$ by virtue of Theorem \ref{thker}, we have {$R_X^\dagger=R_Y^\dagger R_Y\,R_X^\dagger$ so that}
\bea
A_X-A_Y  &  =  &  B\,(R_Y^\dagger S_Y^\tra-R_X^\dagger \,S_X^\tra) =  B\,R_Y^\dagger (S_Y^\tra-R_Y\,R_X^\dagger \,S_X^\tra). \label{sec}
\eea  
Plugging 
\be
\label{eqsy}
S_Y^\tra=B^\tra  Y A\!+\!S^\tra=B^\tra  \Delta A\!+\!S^\tra\!+\!B^\tra X A =B^\tra \Delta A\!+\!S_X^\tra
\ee
and
\be
\label{eqry}
R_Y =R+ B^\tra \,Y\,B-B^\tra\,X\,B+B^\tra\,X\,B=R_X+B^\tra \,\Delta\,B
\ee
into (\ref{sec}) yields
$A_X-A_Y =  B\,R_Y^\dagger (B^\tra\,\Delta\,A-B^\tra\,\Delta\,B\,R_X^\dagger \,S_X^\tra)= B\,R_Y^\dagger B^\tra\,\Delta\,A_X$. This means that $A_Y=A_X-B\,R_Y^\dagger B^\tra\,\Delta\,A_X$.
We already know that in a suitable basis of the state space such that the first coordinates span $\gR^\star_{\ker X}$ and a suitable orthogonal basis of the input space such that the second group of coordinates span $\ker R_X$, matrices $A_X$ and $B$ can be written as in (\ref{dec}), see Lemma \ref{lemmaXB}. The reachable subspace of the pair $(A_X,B\,G_X)$ is written in this basis as $\gR_0=\bsmat I \\[1mm] O \esmat$. 
We want to show that, in the same basis, we also have
\be
\label{dec1}
A_Y=\bmat{cc} A_{X,11} & A_{Y,12} \\ O & A_{Y,22} \emat,\quad \textrm{and} \quad B_2=\bmat{cc} B_{21} \\ O \emat.
\ee
In fact, if this is the case, the reachable subspace of the pair $(A_Y,B\,G_Y)$ exactly coincides with the reachable subspace of the pair $(A_X,B\,G_X)$, i.e., with $\gR_0$, because the pair $(A_{Y,22},0)$ is completely non controllable. In the chosen basis, the difference $\Delta=Y-X$ can be written as
$\Delta=\diag\{O,\Delta_2\}$ in view of \cite[Theorem 2.10]{Stoorvogel-S-98}. \footnote{The result in \cite[Theorem 2.10]{Stoorvogel-S-98} is shown in a basis that is the same considered here. Indeed, the basis of the state space considered in \cite[Theorem 2.10]{Stoorvogel-S-98} has the first coordinates spanning the largest controllability subspace of the quadruple $(A,B,S_X^\tra,R_X)$. However, this subspace coincides with the largest controllability subspace of a quadruple obtained from the previous one by applying the control input $u_t=-K_X\,x_t+H_X\,v_t$, where $\ima H_X=\ker R_X$. The quadruple thus obtained is exactly $(A_X,B\,G_X,S_X^\tra-R_X\,R_X^\tra S_X^\tra,0)=(A_X,B\,G_X,0,0)$, and the corresponding largest controllability subspace is indeed $\gR^\star_{\ker X}$.}
Thus
\beann
A_Y &  =  &  A_X-B\,R_Y^\dagger B^\tra\,\bmat{cc} O & O 	\\ O & \Delta_{2} \emat \bmat{cc} A_{X,11} & A_{X,12} 	\\ O & A_{X,22}\emat \\
 &  =  &  A_X-\bmat{cc} \star \; & \;\star	\\ \star \;& \;\star \emat \bmat{cc} O & O 	\\ O & \Delta_{2}\,A_{X,22}\emat= \bmat{cc} A_{X,11} & \star	\\ O & \star \emat.
\eeann
Therefore, the reachable subspace of the pair $(A_Y,B\,G_Y)$ is exactly $\gR_0 =\bsmat I \\[1mm] O \esmat$, and $A_X|_{\gR^\star_{\ker X}}=A_Y|_{\gR^\star_{\ker Y}}=A_{X,11}$.
\endproof
\ 

We conclude this section by briefly summarizing what we have obtained so far. We have identified a subspace $\gR_0$ such that the closed-loop matrix restricted to this subspace is independent of the particular solution of CGDARE($\Sigma$). This means that if this part of the spectrum contains unstable eigenvalues, CGDARE($\Sigma$) does not admit stabilising solutions for the associated optimal control problem. However, this subspace is a controllability subspace, so it always admits a stabilising friend. This consideration, together with the fact that when $\gR_0$ is non-zero the optimal control is not unique and is parameterised as in (\ref{optcontr}), will lead to the interesting result that the closed-loop can be stabilised exploiting the additional term $G_X$ in (\ref{optcontr}), even in cases in which CGDARE($\Sigma$) does not admit stabilising solutions.

 \section{Stabilisation}
 In the previous sections, we have observed that the eigenvalues of the closed-loop matrix $A_X$ restricted to the subspace $\gR_0$ are independent of the particular solution $X=X^\tra$ of CGDARE($\Sigma$) considered. This means that these eigenvalues -- which, as we will show elsewhere \cite{Ferrante-N-13}, do not appear as generalised eigenvalues of the extended symplectic pencil -- are present in the closed-loop regardless of the solution $X=X^\tra$ of CGDARE($\Sigma$) that we consider.
 On the other hand, we have also observed that $\gR_0$ coincides with the subspace $\gR^\star_{\ker X}$, which is by definition the smallest $(A-B\,K_X)$-invariant subspace containing $\ker X \cap B\,\ker D=\ima (B\,G_{X})$. It follows that it is always possible to find a matrix $L$ that assigns all the eigenvalues of the map $(A_X+B\,G_X\,L)$ restricted to the reachable subspace $\gR^\star_{\ker X}$, by adding a further term $B\,G_X\,L\,x_k$ to the feedback control law, because this does not change the value of the cost with respect to the one obtained by $u_k=-K_X\,x_k$. Indeed, the additional term only affects the part of the trajectory on $\gR^\star_{\ker X}$ which is output-nulling. However, in doing so it may stabilise the closed-loop if $\ker X$ is externally stabilised by $-K_X$.
 We show this fact in the following example.

\begin{example}
{\em 
 Consider a Popov triple in which
$A=\bsmat 1 & 1 \\[1mm] 0 & 1 \esmat$, $B=\bsmat 2 & 0 \\[1mm] 1 & 1 \esmat$, $Q=\bsmat  0 & 0 \\[1mm] 0 & 1  \esmat$, $S=\bsmat 0 & 0 \\[1mm] 0 & 0 \esmat$ and $R=\bsmat  0 & 0 \\[1mm] 0 & 0  \esmat$. The matrix $X=\diag\{0,1\}$ is the only solution of GDARE($\Sigma$) but not a solution of DARE($\Sigma$), since $R+B^\tra\,X\,B$ is singular. Hence, DARE($\Sigma$) does not admit solutions. The corresponding closed-loop matrix $A_X$ is 
$A_X=\diag\{1,0\}$, so that the resulting closed-loop system is not asymptotically stable. However, the solution $X$ of GDARE($\Sigma$) is optimal for the LQ problem, because it leads to the cost $J^\ast=x_2^2(0)$ which cannot be decreased. 
Now, consider the gain $K=B^{-1}A$. This gain leads to the closed-loop matrix $A_{\rm CL}=A-B\,K=0$, and the value of the performance index associated with this closed-loop is again $J=x_2^2(0)=J^\ast$. Therefore, this is another optimal solution of the LQ problem, which differently from $X$ is also stabilising. However, this optimal solution is not associated with any solution of GDARE($\Sigma$), since as aforementioned $X$ is the only solution of GDARE($\Sigma$). In other words, this example shows that there exists an optimal control which is stabilising, but no stabilising solutions of GDARE($\Sigma$) exist. This fact can be explained on the basis of the fact that 
the set of all solutions of the infinite-horizon LQ problem is given by
{\[
\gU_k=\{-K_X\,x_k+G_X v_k \,|\;\; v_k \in\real^m\},
\]
}
where $X$ {is the optimizing solution of GDARE($\Sigma$) and
$G_X=(I_m-R_X^\dagger\,R_X)=\frac{1}{2}\left[ \begin{smallmatrix} 1 & -1 \\[1mm] -1 & 1\end{smallmatrix} \right]$. Therefore, the problem becomes that of using the degree of freedom given by $v_k$} in order to find a closed-loop solution that is optimal and also stabilising. In other words, we determine a matrix $L$ in
{
\be\nonumber
x_{t+1}=(A-B\,K_X)\,x_t+B\,G_X\,L\,x_t=A_X\,x_t+B\,G_X\,L\,x_t
\ee
 such that  the closed-loop
$A_{\rm CL}=A_X+B\,G_X\,L$ is stabilised. 
{It is easy to see that, in general, the set of all optimal closed loop matrices $A_{\rm CL}=A_X+B\,G_X\,L$ are parameterised by}  $A_{\rm CL}=\left[ \begin{smallmatrix} \alpha & \beta \\[1mm] 0 & 0\end{smallmatrix} \right]$ where $\alpha$ and $\beta$ can be arbitrarily chosen by selecting a suitable $L$. In fact, since $B\,G_X=\left[ \begin{smallmatrix} 1 & -1 \\[1mm] 0 & 0\end{smallmatrix} \right]$, by choosing $L=\left[ \begin{smallmatrix} \alpha-1 & \beta \\[1mm] 0 & 0 \end{smallmatrix} \right]$} we obtain the desired form for the closed-loop matrix. Hence, in particular, we can obtain a zero or nilpotent closed-loop matrix. In both cases, the cost is the same and is equal to $J^\ast=x_2^2(0)$.

In other words, there is only one solution to GDARE($\Sigma$) and is not stabilising, and all the optimal solutions of the optimal control problem are given by
 the closed-loop matrix $A_X+B\,G_X\,L$, where $L$ is a degree of freedom. By using this degree of freedom, we have found solutions of the optimal control problem that are stabilising but which do not correspond to stabilising solutions of GDARE($\Sigma$), because GDARE($\Sigma$) does not have stabilising solutions. \endex
}
\end{example}

\section*{Concluding remarks}
In this paper we presented a self-contained analysis of some structural properties of the generalised algebraic Riccati equation that arises in infinite-horizon discrete linear quadratic optimal control. Important side results on Hermitian Stein equations and on spectral factorisation have been established to the end of showing the fundamental role that the term $R_X$ plays in the structure of the solutions of the CGDARE and of the corresponding LQ problem. The considerations that emerged from this analysis have in turn been used to show that a subspace $\gR_0$ can be identified that is independent of the particular solution of CGDARE considered. Even more importantly, it has been shown that the closed-loop matrix restricted to this subspace does not depend on the particular solution of CGDARE. This structural property of GDARE can also be displayed via a decomposition on the extended symplectic pencil, which shows that the spectrum of the closed-loop matrix restricted to $\gR_0$ is not reflected on the generalised eigenstructure of the extended symplectic pencil. In other words, this part of the spectrum has been shown to be fixed for any state-feedback control constructed from a solution of the CGDARE. On the other hand, if such subspace is not zero, i.e., when the related extended symplectic pencil is not regular, in the optimal control  a further term can be added to the state-feedback generated from the solution of the Riccati equation that does not modify the value of the cost. This term can in turn be expressed in state-feedback form, and acts as a degree of freedom that can be employed to stabilise the closed-loop even in cases in which no stabilising solutions exists of the Riccati equation.

The results presented here, | as it will shown in a forthcoming paper \cite{Ferrante-N-13} | can be used to generalise the approach taken in \cite{Ferrante-N-06} to the case of non-regular extended symplectic pencil for the solution of constrained finite-horizon LQ problems.

\end{document}